\documentclass{article}
\usepackage[utf8]{inputenc}
\usepackage{amsmath}
\usepackage{color}
\usepackage[top=2.54cm, bottom=2.54cm, left=2.54 cm, right=2.54 cm]{geometry}
\usepackage{amssymb}
\usepackage{hyperref}
\usepackage[shortlabels]{enumitem}

\usepackage{algorithm,amsthm, forest}
\usepackage{float}
\usepackage{algpseudocode} 
\algnewcommand{\LeftComment}[1]{\Statex \(\triangleright\) #1}

\newtheorem{theorem}{Theorem}[section]
\newtheorem{lemma}[theorem]{Lemma}
\newtheorem{corollary}[theorem]{Corollary}

\newtheorem{definition}[theorem]{Definition}
\newtheorem{conjecture}[theorem]{Conjecture}
\newtheorem{remark}[theorem]{Remark}
\newtheorem{proposition}[theorem]{Proposition}
\newtheorem{question}[theorem]{Question}
\newtheorem{example}[theorem]{Example}

\def\Po{{\textbf{Po}}}
\def\Bin{{\textbf{Bin}}}
\def\G{{\mathcal G}}
\def\I{{\mathcal I}}

\def\pr{{\mathbb P}}
\def\ex{{\mathbb E}}
\def\bfd{{\bf d}}
\def\bft{{\bf t}}

\def\bfg{{\bf g}}

\def\eps{{\epsilon}}
\def\qed{~~\vrule height8pt width4pt depth0pt}

\def\proof{\noindent{\em Proof.~}~}

\newcommand\norm[1]{\|#1\|}
\newcommand\eqn[1]{{(\ref{#1})}}

\newcommand\jt[1]{#1}
\newcommand\fix[1]{#1}

\title{Embedding theorems for random graphs with specified degrees}
\date{}
\author{Pu Gao\thanks{Research supported by NSERC RGPIN-04173-2019.} \\ University of Waterloo \\ pu.gao@uwaterloo.ca \and Yuval Ohapkin \\ University of Waterloo \\ yohapkin@uwaterloo.ca}

\begin{document}

\maketitle

\begin{abstract}
Given an $n\times n$ symmetric matrix $W\in [0,1]^{[n]\times [n]}$, let $\G(n,W)$ be the random graph obtained by independently including each edge $jk\in\binom{[n]}{2}$ with probability $W_{jk}=W_{kj}$. Given a degree sequence $\bfd=(d_1,\ldots, d_n)$, let $\G(n,\bfd)$ denote a uniformly random graph with degree sequence $\bfd$.
We couple $\G(n,W)$ and $\G(n,\bfd)$ together so that asymptotically almost surely $\G(n,W)$ is a subgraph of $\G(n,\bfd)$, where $W$ is some function of $\bfd$. Let $\Delta(\bfd)$ denote the maximum degree in $\bfd$. Our coupling result is optimal when $\Delta(\bfd)^2\ll \norm{\bfd}_1$, i.e.\ $W_{ij}$ is asymptotic to $\pr(ij\in \G(n,\bfd))$ for every $i,j\in [n]$. We also have coupling results for $\bfd$ that are not constrained by the condition $\Delta(\bfd)^2\ll \norm{\bfd}_1$. For such $\bfd$ our coupling result is still close to optimal, in the sense that $W_{ij}$ is asymptotic to $\pr(ij\in \G(n,\bfd))$ for most pairs $ij\in \binom{[n]}{2}$.

\end{abstract}

\section{Introduction}

Given a realisable degree sequence ${\bf d} = (d_1, \ldots, d_n)$ (a degree sequence $\bfd$ is realisable if there exists a simple graph with degree sequence $\bfd$), let $\G(n,\bfd)$ denote a random graph chosen uniformly from the set of graphs on $[n]$ where vertex $i$ has degree $d_i$. 
Random graphs with a specified degree sequence are a popular class of random graphs used in many fields of research \jt{such as social network modelling and analysis~\cite{morris2014random,lusher2013exponential}, epidemic analysis~\cite{bohman2012sir}, and network sciences~\cite{van2024random,chung2002average,newman2001random}.}  While these random graphs have been extensively used to model and analyse real-world networks, such as social networks and the internet, they present several analytical challenges. The most prominent difficulties in analysing $\G(n,\bfd)$ are evaluating the edge probabilities and dealing with edge dependencies. Compared with the classical Erd\H{o}s-R\'{e}nyi random graph $\G(n,p)$ where every edge $jk\in \binom{[n]}{2}$ appears independently with probability $p$, there is no known closed formula for $\pr(jk\in \G(n,\bfd))$ for general $\bfd$. Although asymptotic formulas exist for $\pr(jk\in \G(n,\bfd))$ for some classes of $\bfd$, the correlation between the edges poses additional challenges when estimating the probabilities of events in $\G(n,\bfd)$. 

We denote $\G(n,\bfd)$ by $\G(n,d)$ when $\bfd=(d,d,\ldots,d)$, i.e.\ $\bfd$ is a $d$-regular degree sequence. The random regular graph $\G(n,d)$ is the most well-understood model among all $\G(n,\bfd)$. \jt{For instance, the asymptotic enumeration of $d$-regular graphs on $n$ vertices has been completely solved following a sequence of benchmark research~\cite{bender1978asymptotic,bollobas1980probabilistic,mckay1985asymptotics,mckay1991asymptotic,mckay1990asymptotic,barvinok2010number,barvinok2013number,liebenau2017asymptotic}. Properties and graph parameters of random regular graphs have been extensively studied, including Hamiltonicity, the chromatic and list chromatic numbers, independence number, and the distribution of subgraphs in $\G(n,d)$~\cite{cooper2002random,cooper2002random,krivelevich2001random,gao2008distribution}. We refer the interested readers to a survey of Wormald~\cite{wormald1999models} for many other properties of random regular graphs.}  However, in general $\bfd$ can vary from near-regular sequences to heavy-tailed sequences such as power-law sequences. \jt{Both enumerating graphs of given degree sequences and analysing properties such as connectivity or Hamiltonicity of $\G(n,\bfd)$ turn out to be challenging, and there are many open problems in this field.} In 2004 Kim and Vu proposed the sandwich conjecture, which informally says that $\G(n,d)$ can be well approximated by $\G(n,p=d/n)$ through sandwiching $\G(n,d)$ between two correlated copies of $\G(n,p)$, one with $p$ slightly smaller than $d/n$ and the other with $p$ slightly greater than $d/n$. To formally state the conjecture, we define a coupling of a finite set of random variables (or graphs) $Z_1,\ldots, Z_k$ as a $k$-tuple of random variables $(\hat Z_1,\ldots,\hat Z_k)$ defined in the same probability space such that the marginal distribution of $\hat Z_i$ is the same as the distribution of $Z_i$ for every $1\le i\le k$.
The sandwich conjecture \jt{is about a 3-tuple coupling of random graphs, proposed} by Kim and Vu,  given as below. \jt{The standard Landau notation is used, and a formal definition of the notation is given before Section~\ref{sec:CL}.}
\begin{conjecture}[\cite{kim2004sandwiching}]
For $d\gg \log n$, there exist probabilities $p_1,p_2=(1+o(1))d/n$ and a coupling $(G_L, G, G_U)$ such that marginally, $G_L\sim \G(n,p_1)$, $G\sim \G(n,d)$, $G_U\sim \G(n,p_2)$ and jointly, $\pr(G_L\subseteq G\subseteq G_U)=1-o(1)$.
\end{conjecture}

The sandwich conjecture, if proved to be true, is a powerful tool for analysing $\G(n,d)$, and reveals beautiful distributional relations between the two different random graph models. The assumption $d\gg \log n$ is necessary, as for $p=O(\log n/n)$, the maximum and minimum degrees of $\G(n,p)$ differ by some constant factor away from 1, making it impossible for the sandwich conjecture to hold. For simplicity, we refer to a coupling $(\G_n^1,\G_n^2)$ of two random graphs $\G_n^{1}$ and $\G_n^{2}$ as {\em embedding $\G_n^1$ into $\G_n^2$}, if in the coupling  $\G_n^1$ is a subgraph of $\G_n^2$ with probability going to 1 as $n\to\infty$. It turns out that embedding $\G(n,d)$ into $\G(n,p=(1+o(1))d/n)$ is much more difficult than embedding $\G(n,p=(1-o(1))d/n)$ into $\G(n,d)$. In their paper~\cite{kim2004sandwiching}, Kim and Vu established an embedding of $\G(n,p=(1-o(1))d/n)$ into $\G(n,d)$ for $d$ up to approximately $n^{1/3}$ (and $d\gg \log n$).  Subsequent work by Dudek, Frieze, Ruci\'{n}ski and  {\v{S}ileikis~\cite{dudek2017embedding} improved this result to $d=o(n)$, and extended it to random uniform hypergraphs as well. The first 2-side sandwich theorem was proved by Isaev, McKay and the first author~\cite{gao2020sandwiching,gao2022sandwiching} in the case that $d$ is linear in $n$. Klimo\"{s}ov\'{a}, Reiher, Ruci\'{n}ski and  {\v{S}ileikis~\cite{klimovsova2023sandwiching} proved the sandwich conjecture for $d\gg (n\log n)^{3/4}$. \jt{Finally Isaev, McKay and the first author~\cite{gao2020kim} confirmed the sandwich conjecture for all $d\ge \log^4 n$.} It is worth noting that a more general sandwich theorem was presented in~\cite{gao2020sandwiching,gao2022sandwiching,gao2020kim} for random graphs $\G(n,\bfd)$ where $\bfd$ is a near-regular degree sequence, where all degrees are asymptotically equal. 

While it is not possible to approximate $\G(n,\bfd)$ by $\G(n,p)$ in a useful way for more general forms of $\bfd$ \jt{(since the typical degree distribution of $\G(n,p)$ would be very different from $\bfd$),} a natural alternative is to consider a generalised Erd\"{o}s-R\'{e}nyi graph $\G(n,W)$, where $W$ is a symmetric $n\times n$ matrix. In this model, each edge $jk\in\binom{[n]}{2}$ appears in $\G(n,W)$ independently with probability $W_{jk}$. By choosing different forms of $W$, we can recover well-studied models in the literature, such as $\G(n,p)$, the Chung-Lu model~\cite{chung2002connected,chung2002average,chung2004spectra,chung2006volume}, the $\boldsymbol{\beta}$-model~\cite{chatterjee2011random,isaev2018complex}, and the stochastic block model~\cite{holland1983stochastic}.

The objective of this paper is to establish an embedding of $\G(n,W)$ into $\G(n,\bfd)$, where $\bfd$ is  a degree sequence that may deviate significantly from being regular, for a suitable choice of $W$. What constitutes a good $W$? Intuitively one would like the entry $W_{jk}$ to be approximately $\pr(jk\in \G(n,\bfd))$, the true probability that $jk$ is an edge in $\G(n,\bfd)$. To formalise this notion, we define $W^*=W^*(\bfd)$ as the $n\times n$ matrix where $W^*_{jk}=\pr(jk\in \G(n,\bfd))$ for every $i,j\in [n]$. 

\begin{question}\label{question}
Given $\bfd$, does there exists $W=(1-o(1))W^*(\bfd)$ such that $\G(n,W)$ can be embedded into $\G(n,\bfd)$?
\end{question}

 Note that the asymptotic value of $\pr(jk\in \G(n,\bfd))$ is not known for all realisable degree sequences $\bfd$. The coupling scheme we use in this paper heavily relies on the expression of the asymptotic formula for $\pr(jk\in \G(n,\bfd))$. Therefore, we are restricted to degree sequences for which an asymptotic evaluation of $\pr(jk\in \G(n,\bfd))$ is available in the literature.  
Without loss of generality we may assume that $d_1 \ge \cdots \ge d_n \ge  1$. Let $\Delta(\bfd) = d_1$ and $\delta(\bfd) = d_n$. Another important parameter of $\bfd$ is $J(\bfd)=\sum_{i=1}^{d_1} d_i$. The following \fix{proposition, a direct corollary of~\cite[Theorem 1]{gao2020subgraph},
 estimates the edge probabilities in $\G(n,\bfd)$ when $J(\bfd)=o(\norm{\bfd}_1)$.


\begin{proposition} (Corollary of~\cite[Theorem 1]{gao2020subgraph}) \label{p:edgeprob}
Suppose that $J(\bfd)=o(\norm{\bf d}_1)$. Then
\[
\pr(jk \in \G(n, {\bf d}) ) =\left(1+O\left(\frac{J(\bfd)}{\norm{\bf d}_1}\right)\right) \frac{d_jd_k}{\norm{\bfd}_1+d_j d_k}.
\]
\end{proposition}
}

\jt{
\begin{remark}
The condition $J(\bfd)=o(\norm{\bfd}_1)$ requires $\bfd$ to be a sparse degree sequence. It also restrict the tail of $\bfd$ so that there cannot be too many degrees that are very large. See Proposition~\ref{p:properties} in the Appendix for properties of $\bfd$ that satisfies this technical condition. The edge probability in Proposition~\ref{p:edgeprob} is in general not correct without the condition $J(\bfd)=o(\norm{\bfd}_1)$. See an example like this in Remark~\ref{r:edgeprob}.
\end{remark}

}

\begin{remark} \label{r:examples}
\fix{Although Proposition~\ref{p:edgeprob} does not give asymptotic estimate of $\pr(jk \in \G(n, {\bf d}) )$ for every degree sequence,}
there are many interesting families of degree sequences satisfying $J(\bfd)=o(\norm{\bfd}_1)$, including examples of 
\begin{itemize}
\item all $d$-regular, and near-$d$-regular degree sequences \jt{(namely, all degrees are asymptotic to $d$)} where $d=o(n)$;
\item all perturbed sequences from a near-$d$-regular degree sequence by arbitrarily decreasing at most $(1-c)n$ entries, where $c>0$ is fixed, and $d=o(n)$; 
\item heavy-tailed degree sequences such as power-law sequences. We refer interested readers to~\cite[Section 2]{gao2016enumeration} for a \jt{formal definition of power-law sequences and other examples of heavy-tailed degree sequences.}
\end{itemize}

\fix{Assume that $J(\bfd)=o(\norm{\bfd}_1)$. It is easy to show that for most pairs $jk$, $d_jd_k=o(\norm{\bfd}_1)$ (e.g.\ see Proposition~\ref{p:properties} in the Appendix), and thus  by Proposition~\ref{p:edgeprob} their edge probabilities are $o(1)$.} On the other hand, if there exist pairs $jk$ such that $d_jd_k=\Omega(\norm{\bfd}_1)$, then Proposition~\ref{p:edgeprob} implies that the edge probabilities in $\G(n,\bfd)$ \fix{will not be uniformly $o(1)$. In particular, for every such pair $jk$ the edge probability of $jk$ is bounded from below by some constant $c>1$. In fact, there exist degree sequences $\bfd$ satisfying $J(\bfd)=o(\norm{\bfd}_1)$ whose edge probabilities} takes values ranging from $o(1)$ to constant $0<c<1$  and to $1-o(1)$. \jt{The power-law sequences mentioned in Remark~\ref{r:examples} are good examples like this.}  As we will see later, it is the presence of such diverse edge probabilities  that poses a challenge to embedding $\G(n,W)$ into $\G(n,\bfd)$.
\end{remark}

In light of Proposition~\ref{p:edgeprob}, we define matrix $P(\bfd)$, which is asymptotic to $W^*(\bfd)$ under the condition $J(\bfd)=o(\norm{\bfd}_1)$.
\begin{definition}\label{def:P}
Given ${\bf d}=(d_1,\ldots, d_n)$, let $P({\bfd})$ be the symmetric $n\times n$ matrix defined by $P_{ij}=P_{ji}=\frac{d_id_j}{\|\bfd\|_1+d_id_j}$, for every $1\le i<j\le n$, and $P_{ii}=0$ for every $1\le i\le n$.
\end{definition}
One of our main results is the following theorem, which gives a positive answer to Question~\ref{question} for degree sequences satisfying $\Delta({\bf d})^2 = o(\norm{\bfd}_1)$, a condition that is stronger than $J(\bfd)=o(\norm{\bfd}_1)$.

\begin{theorem}\label{thm:main}
Assume that ${\bf d} = {\bf d}(n)$ is a degree sequence satisfying $\Delta({\bf d})^2 = o(\norm{\bfd}_1)$ and $\delta({\bf d}) \gg\log{n}$. Then, there exist \jt{$\varepsilon=o(1)$} and a coupling $(G_L, G)$, where $G_L \sim \G(n, W)$, \jt{$W=(1-\varepsilon)P(\bfd)$,} and $G \sim \G(n, {\bf d})$, such that $\pr(G_L \subseteq G) = 1-o(1)$.
\end{theorem}

Under the stronger condition that $\Delta(\bfd)^2=o(\norm{\bfd}_1)$, \jt{$d_jd_k=o(\norm{\bfd}_1)$ for all $jk\in \binom{[n]}{2}$, and thus} $\pr(jk\in \G(n,\bfd))=o(1)$ for all $jk\in \binom{[n]}{2}$ by Proposition~\ref{p:edgeprob}. This property is a crucial element in our coupling scheme which allows us to achieve an optimal embedding. Our next theorem embeds $\G(n,W)$ into $\G(n,\bfd)$ for $\bfd$ satisfying $J(\bfd)=o(\norm{\bfd}_1)$. In this case, $W_{jk}=(1-o(1))P(\bfd)_{jk}$ for all $jk$ where $\pr(jk\in \G(n,\bfd))=o(1)$. However, $W_{jk}$ is not asymptotic to $\pr(jk\in \G(n,\bfd))$ for $jk$ where $\pr(jk\in \G(n,\bfd))$ is bounded away from 0.  Given a function $f: {\mathbb R}\to {\mathbb R}$, let $f(P)$ denote the matrix $(f(P_{i,j}))_{i,j\in [n]}$.

\begin{theorem}\label{thm:main2}
Assume that ${\bf d} $ is a degree sequence satisfying $J(\bfd)=o(\norm{\bfd}_1)$ and $\delta({\bf d}) \gg\log{n}$. Then, there exist \jt{$\varepsilon=o(1)$} and a coupling $(G_L, G)$, where $G_L \sim \G(n, W)$, \jt{$W=(1-\varepsilon)f(P(\bfd))$} where $f(x)=1-e^{-x}$, and $G \sim \G(n, {\bf d})$, such that $\pr(G_L \subseteq G) = 1-o(1)$.
\end{theorem}

We see that Theorem~\ref{thm:main} immediately follows from Theorem~\ref{thm:main2}. \smallskip

\noindent {\em Proof of Theorem~\ref{thm:main}.} \jt{Note that $\Delta(\bfd)^2=o(\norm{\bfd}_1)$ implies $J(\bfd)=o(\norm{\bfd}_1)$ since $J(\bfd)\le \Delta(\bfd)^2$. Moreover, $\Delta(\bfd)^2=o(\norm{\bfd}_1)$ implies that $P({\bfd})_{jk}=o(1)$ for all $jk$ by Definition~\ref{def:P}.} Now Theorem~\ref{thm:main} follows as a straightforward corollary of Theorem~\ref{thm:main2} by noticing that $1-e^{-x}=(1+O(x))x$. \qed \medskip

\begin{remark}
\label{r:edgeprob}
Note that $P(\bfd)=(1+o(1))W^*(\bfd)$ is not true in general. For instance, consider the $d$-regular degree sequence where $d=\Theta(n)$. \jt{By symmetry we know that $W^*_{jk}=d/(n-1)\sim d/n$ for every $jk$. On the other hand, by Definition~\ref{def:P}, $P_{jk}=d^2/(dn+d^2)=d/(n+d)$ which is not asymptotic to $d/n$ as $d=\Theta(n)$. Thus, in this example, $P_{jk}$ is significantly smaller than $W^*_{jk}$.} There are two open problems \jt{below} that can be very interesting for future research. \jt{The first question is a weaker version of Question~\ref{question}, replacing $W^*(\bfd)$ by $P(\bfd)$.}
\begin{enumerate}
\item[(a)] Does Theorem~\ref{thm:main2} hold with some $W=(1+o(1))P(\bfd)$?
\end{enumerate}
\jt{Note that Theorem~\ref{thm:main} answers the question above only for degree sequences satisfying $\Delta^2=o(\norm{\bfd}_1)$. However, in the more general setting  as Theorem~\ref{thm:main2}, there can be pairs $jk$ where $d_jd_k=\Omega(\norm{\bfd}_1)$, and for such pairs, our current embedding requires $W$ where $W_{jk}$ is significantly smaller than $P_{jk}$. Hence we do not have a complete answer to question (a). }
\begin{enumerate}
\item[(b)] Can Question~\ref{question} be answered for $\bfd$ \jt{where $J(\bfd)=o(\norm{\bfd}_1)$ is not satisfied}? In particular, is there a way to embed $\G(n,W)$ into $\G(n,\bfd)$ for some reasonable $W$ without knowing the asymptotic value of $W^*(\bfd)$ or the conditional edge probabilities as in Proposition~\ref{p:edgeprob}?
\end{enumerate}

\end{remark}

\begin{remark}
To obtain a sandwich type result as for $\G(n,d)$, we would hope to embed $\G(n,\bfd)$ into $\G(n, W)$ for some $W=(1+o(1))W^*$. If we approach in the same way as for $\G(n,d)$, we would embed $\G(n, (J-I)-(1+o(1))W^*)$ into $\G(n,(n-1){\bf 1}-\bfd)$ where ${\bf 1}$ is the all-one vector, $J={\bf 1}{\bf 1}^T$ and $I$ is the identity matrix. Both the proofs in~\cite{gao2020sandwiching,klimovsova2023sandwiching} for embedding $\G(n,1-(1+o(1))d/n)$ into $\G(n,n-1-d)$ use some counting arguments that  heavily  rely on the fact that the underlying graphs (during the construction of the coupling) have almost equal degrees, \jt{and we do not think this argument can extend} to graphs that are far away from being regular. In fact, we do not have enough intuition to support a  sandwich conjecture for $\G(n,\bfd)$ as for $\G(n,d)$, especially for degree sequences $\bfd$ where the values of edge probabilities range from $o(1)$ to $1-o(1)$.   
\end{remark}

Throughout the paper $n$ is assumed to be sufficiently large. 
For two sequences of real numbers $(a_n)_{n=0}^{\infty}$ and $(b_n)_{n=0}^{\infty}$, we write $a_n=O(b_n)$ if there is $C>0$ such that $|a_n|<C|b_n|$ for every $n$. We say $a_n=\Omega(b_n)$ if $a_n>0$ and $b_n=O(a_n)$. We say $b_n = o(a_n)$, or $b_n \ll a_n$, or $a_n \gg b_n$ if $a_n>0$ and $\lim_{n \rightarrow \infty} b_n/a_n = 0$. We say a sequence of events $A_n$ indexed by $n$ holds a.a.s.\ (asymptotically almost surely) if $\pr(A_n)=1-o(1)$. All asymptotics are with respect to $n \rightarrow \infty$.

\subsection{Relation between $\G(n,P(\bfd))$ and the Chung-Lu model}
\label{sec:CL}
In 2002, Chung and Lu introduced the random graph model with given expected degrees~\cite{chung2002connected,chung2002average,chung2004spectra,chung2006volume}. Given a sequence of nonnegative real numbers ${\bf w}=(w_1,\ldots, w_n)$ satisfying that 

\begin{equation}
\max_i w_i^2 \le \norm{{\bf w}}_1, \label{assumption}
\end{equation}
the random graph $G({\bf w})$ is defined by $\G(n,\hat W({\bf w}))$ where $\hat W_{jk}({\bf w})=w_jw_k/\norm{{\bf w}}_1$. To avoid relying on the assumption~\eqn{assumption} one can define $\hat W_{jk}=\min\{w_jw_k/\norm{{\bf w}}_1,1\}$. However, most of the work about the Chung-Lu model assumed~\eqn{assumption}. Notice that the three matrices $P(\bfd)$, $W^*(\bfd)$ and $\hat W(\bfd)$ are all asymptotically equal if $\bfd$ satisfies $\Delta(\bfd)^2=o(\norm{\bfd}_1)$, \jt{as under this assumption, $P(\bfd)_{jk}\sim W^*(\bfd)_{jk}\sim \hat W(\bfd)_{jk}=d_jd_k/\norm{\bfd}_1$ for all $jk$, by definition of these three matrices and Proposition~\ref{p:edgeprob}}. However, if there exists $jk$ such that $d_jd_k=\Omega(\norm{\bfd}_1)$ then $P(\bfd)_{jk}$, $W^*(\bfd)_{jk}$ and $\hat W(\bfd)_{jk}$ may be all asymptotically distinct.

\subsection{Applications}

Given two $n\times n$ matrices $W_1$ and $W_2$,
we say $W_1\le W_2$ if $(W_1)_{ij}\le (W_2)_{ij}$ for every $i,j\in[n]$. It is well known that $\G(n,W)$ has the nice ``nesting property'', meaning that $\G(n,W_1)$ can be embedded into $\G(n,W_2)$ provided that $W_1\le W_2$. However, $\G(n,\bfd)$ does not have the nesting property. Given two degree sequences $\bfd\preceq \bfg$, it is in general not true that $\G(n,\bfd)$ can be embedded into $\G(n,\bfg)$. In fact, it is easy to construct degree sequences $\bfd\preceq \bfg$ for which there exists $jk$ such that $\pr(jk\in \G(n,\bfd))= 1-o(1)$ and $\pr(jk\in \G(n,\bfg))=o(1)$ \jt{(see an example like this in the Appendix)}.  Consequently and perhaps rather surprisingly, it is difficult to prove some rather ``intuitive'' results about $\G(n,\bfd)$. For instance, although it is known that a.a.s.\ $\G(n,3)$ is connected, to our knowledge it is not known if $\G(n,\bfd)$ is a.a.s.\ connected for every $\bfd$ provided that $\delta(\bfd)\ge 3$~\cite{gao2023subgraph}. Most such results are restricted to certain families of degree sequences for which either some enumeration results are known, or some enumeration proof techniques can be applied. Theorem~\ref{thm:main2} gives a powerful tool to obtain such results by first embedding $\G(n,W)$ into $\G(n,\bfd)$ and then applying the nesting property of $\G(n,W)$. We show a few examples below.

\begin{theorem}
Assume ${\bf d} $ is a degree sequence satisfying $J(\bfd)=o(\norm{\bfd}_1)$ and $\delta({\bf d}) \gg\log{n}$. Then,  
\begin{enumerate}[(a)]
\item if $\delta(\bfd)^2/\norm{\bfd}_1\ge (1+c)\log n/n$ for some fixed $c>0$ then a.a.s.\
 $\G(n,\bfd)$ is Hamiltonian and $k$-connected for every fixed $k$;
\item if $\delta(\bfd)^2/\norm{\bfd}_1\ge c/n$ for some fixed $c>1$ then a.a.s.\ $\G(n,\bfd)$ contains $H$ as a minor for every fixed graph $H$;
\item if $\delta(\bfd)^2/\norm{\bfd}_1\ge n^{-2/(k+1)+c}$ for some fixed $c>0$ where $k>0$ is a fixed integer, then a.a.s.\ simultaneously for all graphs $H$ on $[n]$ with maximum degree at most $k$, $\G(n,\bfd)$ has a subgraph isomorphic to $H$. 
\end{enumerate}
\end{theorem}

\proof By definition, $P(\bfd)_{jk}\ge (1-o(1))\delta(\bfd)^2/\norm{\bfd}_1$ for every $jk$. Thus, by Theorem~\ref{thm:main2} and the nesting property of $\G(n,W)$, $\G(n,p)$ can be embedded into $\G(n,\bfd)$ for some $p= (1-o(1))\delta(\bfd)^2/\norm{\bfd}_1$. \jt{Parts (a,b) follow as for every fixed $c>0$, $\G(n,(1+c)\log n/n)$ is a.a.s.\ Hamiltonian and $k$-connected~\cite{bollobas1998random}; and $\G(n,c/n)$ contains every fixed graph minor when $c>1$~\cite{fountoulakis2008order}. Part (c) follows from Theorem 7.2 of~\cite{frankston2021thresholds}.}\qed

\begin{remark}
It should be \jt{possible} to improve or even remove some assumptions such as  $\delta(\bfd)^2/\norm{\bfd}_1\ge (1+c)\log n/n$ for some fixed $c>0$ in part (a), by directly analysing $\G(n, (1+o(1))f(P(\bfd)))$. \jt{For instance, the sharp threshold of Hamiltonicity was studied for the stochastic block model~\cite{anastos2021hamiltonicity}.} We did not attempt it as the main objective of this paper is to prove the embedding theorems.
\end{remark}

\section{The coupling procedure}

\subsection{The old and the new}
\label{sec:new}

Assume that we aim to embed $\G(n,p)$ into $\G(n,d)$ where $p=(1-o(1))d/n$.
Regardless of a few minor differences, the coupling procedures employed in~\cite{kim2004sandwiching,dudek2017embedding,gao2020sandwiching,gao2022sandwiching,klimovsova2023sandwiching,gao2020kim} \jt{all use the following approach, which was introduced in the original paper of Kim and Vu~\cite{kim2004sandwiching}}: let $x_1,x_2,\ldots, x_m$ be a sequence of random edges, each uniformly and independently chosen from $K_n$ \jt{(here $m$ is a carefully chosen integer-valued random variable that is concentrated around $\norm{\bfd}_1/2$)}. Sequentially add edges in this sequence to $G$ and to $G_L$ respectively. With a small probability $\eps_i=o(1)$, edge $x_i$ is rejected by $G$; and with a slightly larger \jt{(than $\eps_i$)} but still rather small probability $\zeta$, $x_i$ is rejected by $G_L$. The parameter $\eps_i$ is chosen to be proportional to $p_i(x_i)$, the probability that $x_i$ is an edge of $\G(n,\bfd)$ conditional on the event that all the edges that have been added to $G$ are edges of $\G(n,\bfd)$. The key idea of the coupling procedure is that, until $m$ gets very close to $dn/2$, with high probability, $p_i(jk)$ is approximately the same for every edge $jk$ that has not been added to $G$ yet. Since $x_i$ is uniformly chosen, a small rejection probability $\eps_i$ suffices to ensure that $x_i$ is added to $G$ according to the correct conditional probability. We can prove that with high probability, $\eps_i=o(1)$ for every $1\le i\le m$, where $m$, \jt{as commented earlier, is concentrated around} $(1-o(1))dn/2$. Hence we can choose some $\zeta=o(1)$ such that $\eps_i\le \zeta$ for every $1\le i\le m$. Moreover, (a) $G_L$ obtained after the $m$-th iteration is a uniformly random graph conditional on the number of edges it contains, as every edge in $\binom{[n]}{2}$ has an equal probability to be added to $G_L$; and (b) $G_L$ is a subgraph of $G$, as the rejection probability $\zeta$ for $G_L$ is slightly larger than that for $G$ in every step.

Now we are considering $\bfd$ where the degrees are not all asymptotically the same. The most natural way to extend the previous coupling procedure is to generate the sequence of i.i.d.\ random edges $x_1,x_2,\ldots, x_m$ where each edge is chosen with probability proportional to $W^*(\bfd)$. This coupling procedure works well when $\Delta(\bfd)^2=o(\norm{\bfd}_1)$. This is because the conditional probability $p_i(jk)$ of $jk$ being an edge of $\G(n,\bfd)$ in step $i$ of the procedure turns out to be  proportional to $W^*(\bfd)_{jk}$ uniformly for all $jk$ during the whole coupling process, resulting a small rejection probability $\eps_i$. However, rather surprisingly, if there exists $jk$ such that $d_jd_k=\Omega(\norm{\bfd}_1)$ then the ratio $p_i(jk)/W^*_{jk}$ changes in a non-uniform way over $jk$ and over time $i$ of the coupling procedure, resulting larger and larger rejection probability $\eps_i$. As $\zeta$ has to be chosen uniformly during the process, we can only choose $\zeta=1-o(1)$, meaning that almost all edges are rejected by $G_L$. The coupling procedure thus fails. \smallskip

To overcome the challenges, we develop two novel ideas.

\begin{itemize}
\item Instead of using the same rejection probability $\zeta$ for every edge to be added to $G_L$, we use different rejection probabilities for different edges. However, for a  given edge $jk$, the rejection probability $\zeta_{jk}$ remains uniform throughout the procedure. This uniformity is needed for the output $G_L$ to have the correct distribution.

\item During the coupling procedure, the value of $p_i(jk)/W^*_{jk}$ decreases significantly for edges $jk$ where $d_jd_k=O(\norm{\bfd}_1)$, but changes little for $jk$ where $d_jd_k\gg \norm{\bfd}_1$. Consequently, many rejections (particularly those where $d_jd_k\ll \norm{\bfd}_1$) occur already for the construction of $G$, and thus the first idea above would not help, as  even more rejections would occur for $G_L$ than for $G$. To reduce the rejection probability for the construction of $G$, we ``intentionally'' boost the probability (by an $\omega(1)$ factor) of generating edges $jk$ where $d_jd_k\gg \norm{\bfd}_1$ in the sequence  of random edges $x_1,x_2,\ldots$. This probability boosting strategy magically reduces the rejection probability $\eps_i(jk)$ for $jk$ such that $d_jd_k\ll \norm{\bfd}_1$ (note that most $jk\in\binom{[n]}{2}$ are this type; \jt{see Proposition~\ref{p:properties} in the Appendix}). However, the rejection probabilities $\eps_i(jk)$ will be high (close to 1) for $jk$ where $d_jd_k \gg\norm{\bfd}_1$. Nonetheless, the first idea mentioned earlier will be applicable in this case --- we reject these $jk$ more often than the others for the construction of $G_L$. Consequently, the edge probability for such $jk$ in the final construction of $G_L$ turns out to be close to $1-e^{-1}$ instead of $1$. \jt{Note that for pairs $jk$ where $d_jd_k \gg\norm{\bfd}_1$, their edge probability $W^*_{jk}$ in $\G(n,\bfd)$ is asymptotic to 1 by Proposition~\ref{p:edgeprob}. Thus, their edge probabilities in $G_L$ (which is asymptotic to $1-e^{-1}$) is smaller than $W^*_{jk}$, implying that our embedding is not ``tight''. The good news is that}  only a small proportion of $jk\in\binom{[n]}{2}$ are of this type. 
\end{itemize}

We define the procedure for sequential generation of $\G(n,W)$ in Section~\ref{sec:W}, the procedure for sequential generation of $\G(n,\bfd)$ in Section~\ref{sec:d}. In Section~\ref{sec:coupling} we couple the two procedures together to sequentially generate the coupled pair $(G_L,G)$. The parameters used by the coupling procedure will be specified in Section~\ref{sec:parameters}. Finally the proof of Theorem~\ref{thm:main2} is given in Section~\ref{sec:proof}.


\subsection{Sequential generation of $\G(n,W)$}
\label{sec:W}

We define a sequential sampling procedure \textsc{SeqApprox-P}($\bfd,\lambda,\Lambda$), which adds a sequence of edges one at a time, and outputs a random graph where every edge in $\binom{[n]}{2}$ appears independently. The input $\lambda$ is a positive real number, and $\Lambda\in [0,1]^{[n]\times [n]}$ is a symmetric $n\times n$ matrix.   We may think of $\bfd$ as the target degree sequence. As mentioned before, instead of weighting edge probabilities according to their true probabilities in $\G(n,\bfd)$, we boost the probability of edge $jk$ if $d_jd_k$ is large.  To formalise this notion, we define matrix $Q$ as follows, which provides the probability distribution for the edges to be sequentially sampled in the procedure \textsc{SeqApprox-P}($\bfd,\lambda,\Lambda$).

\begin{definition}
Given ${\bf d}=(d_1,\ldots, d_n)$, let $Q({\bfd})$ be the symmetric $n\times n$ matrix defined by $Q_{ij}=Q_{ji}=(\sum_{1\le k<\ell\le n}d_kd_\ell)^{-1}d_id_j$, for every $1\le i<j\le n$, and $Q_{ii}=0$ for every $1\le i\le n$.
\end{definition}
 
Procedure \textsc{SeqApprox-P}($\bfd,\lambda,\Lambda$) is given below. \jt{We use $\Po(\lambda)$ to denote the Poisson distribution with mean $\lambda$.}

\begin{algorithm}[H]
  \label{alg:W}
  \begin{algorithmic}[1]
    \Procedure{SeqApprox-P}{$\bfd,\lambda,\Lambda$}
    \State Let  $\I \sim \Po(\lambda)$. 
    \State Let $G_0$ be the empty graph on $[n]$. 
     \For{$i$ in $1, \dots,\I$}
     \State Pick an edge $jk\in \binom{[n]}{2}$ with probability proportional to $d_jd_k$ (i.e.\ with probability $Q(\bfd)_{jk}$).
     \State $G_i=G_{i-1}\cup \{jk\}$ with probability $\Lambda_{jk}$, and $G_i=G_{i-1}$ with probability $1-\Lambda_{jk}$. 
     \EndFor
    \State Return $G_{\I}$. 
    \EndProcedure
    \end{algorithmic}
    \end{algorithm}

\jt{The parameters $\lambda$ and $\Lambda$ will be set in Section~\ref{sec:parameters}. Roughly speaking, $\lambda$ is approximately $\norm{\bfd}_1/2$ so that the number of edges in the output of \textsc{SeqApprox-P}($\bfd,\lambda,\Lambda$) is close to the number of edges in $\G(n,\bfd)$, if only a small proportion of them are to be rejected, which we will prove later. The matrix $\Lambda$ is chosen so that $\Lambda_{jk}$ is approximately $\norm{\bfd}_1/(d_jd_k+\norm{\bfd}_1)$. Since each edge $jk$ is selected with probability proportional to $d_jd_k$ by the definition of $Q$, and due to the choice of $\Lambda$, edge $jk$ is selected and accepted with probability proportional to $d_jd_k/(d_jd_k+\norm{\bfd}_1)$, which is approximately $W^*_{jk}$ as desired. The exact values of $\lambda$ and $\Lambda$ are only important when we couple the two sampling processes together. At this point, regardless of their values,} we  prove that the output of  \textsc{SeqApprox-P}($\bfd,\lambda,\Lambda$) has distribution $\G(n,W)$ for some $W$ as a function of $\bfd,\lambda$ and $\Lambda$. For two matrices $A$ and $B$ of the same dimension, we denote by $A\odot B$ the Hadamard product of $A$ and $B$, defined by $(A\odot B)_{ij}= A_{ij} B_{ij}$ for every entry $ij$.

\begin{lemma}\label{lem:P1}
 \textsc{SeqApprox-P}($\bfd,\lambda,\Lambda$) returns a random graph $G$ with distribution $\G(n,f(\Lambda\odot Q))$, where $Q=Q(\bfd)$ and $f(x)=1-\exp(-\lambda x)$.
\end{lemma}

\proof Let $e_1 = j_1k_1, \dots, e_N = j_Nk_N$ be an enumeration of the edges in $\binom{[n]}{2}$, where $N = \binom{n}{2}$.For $1 \le i \le N$, let $X_i$ denote the number of times that edge $e_i$ is sampled throughout \textsc{SeqApprox-P}($\bfd,\lambda,\Lambda$).   \jt{We prove that the components of ${\bf X}=(X_1,\ldots, X_N)$ are mutually independent, and consequently, every edge $e_i$ appears  independently in $G$ with probability $\pr(X_i\ge 1)$.} 
For each edge $e_i\in\binom{[n]}{2}$, the probability that $e_i$ is in $G$ is thus given by
\begin{align*}
\pr(X_i \ge 1) &= 1-\pr(X_i=0) =1-\sum_{m=0}^{\infty} e^{-\lambda} \frac{\lambda^m}{m!} \sum_{j=0}^{m} \binom{m}{j} Q_{e_i}^j (1-Q_{e_i})^{m-j} \left(1-\Lambda_{e_i}\right)^j \\
& =1- \sum_{m=0}^{\infty} e^{-\lambda} \frac{\lambda^m}{m!} (1- Q_{e_i}+Q_{e_i}(1-\Lambda_{e_i}))^m \\
&= 1- \exp\left(-\lambda+\lambda(1- Q_{e_i}+Q_{e_i}(1-\Lambda_{e_i}))\right)= 1-\exp(-\lambda Q_{e_i}\Lambda_{e_i}), 
\end{align*}
and the probability generating function for $X_i$ is  given by
\begin{align*}
\ex z^{X_i} &=\sum_{m=0}^{\infty} e^{-\lambda} \frac{\lambda^m}{m!}\sum_{j=0}^{m} \binom{m}{j} Q_{e_i}^j (1-Q_{e_i})^{m-j}  \sum_{k=0}^j \binom{j}{i} \Lambda_{e_i}^k (1-\Lambda_{e_i})^{j-k} z^k \\
&=\sum_{m=0}^{\infty} e^{-\lambda} \frac{\lambda^m}{m!} \sum_{j=0}^{m} \binom{m}{j} Q_{e_i}^j (1-Q_{e_i})^{m-j} (1-\Lambda_{e_i}+\Lambda_{e_i}z)^j\\
& =\sum_{m=0}^{\infty} e^{-\lambda} \frac{\lambda^m}{m!} (1-Q_{e_i}+Q_{e_i}(1-\Lambda_{e_i}+\Lambda_{e_i}z))^m=\exp(-\lambda Q_{e_i}\Lambda_{e_i}(1-z)).
\end{align*}
On the other hand, the probability generating function for the random vector ${\bf X}$ is given by
\begin{align*}
&\sum_{j_1, \dots, j_N} \pr(X_1 = j_1, \dots, X_N = j_N) z_1^{j_1}\cdots z_N^{j_N} = \sum_{m = 0}^{\infty} e^{-\lambda} \frac{\lambda^m}{m!} \left( \sum\limits_{i=1}^{N} \big( Q_{e_i} \Lambda_{e_i} z_i +Q_{e_i}(1-\Lambda_{e_i})\big) \right)^m\\
&\hspace{1cm}= e^{-\lambda} \exp\left(\sum_{i=1}^N  \big(\lambda Q_{e_i}\Lambda_{e_i}z_i +\lambda Q_{e_i}(1-\Lambda_{e_i}) \big)\right) =  \prod_{i=1}^N \exp\left(-\lambda \ Q_{e_i} \Lambda_{e_i}+\lambda Q_{e_i} \Lambda_{e_i} z_i\right)=  \prod_{i=1}^N  \ex z^{X_i},
\end{align*}
where the second last equation above holds because $\sum_{i=1}^N Q_{e_i}=1$ by definition of $Q$.
Hence we have shown that the components of ${\bf X}$ are mutually independent. Thus, $G\sim \G(n,W)$ where 
$W$ is the symmetric $n\times n$ matrix given by $W_{ij}=1-\exp(-\lambda Q_{ij}\Lambda_{ij})$ for every $i,j\in [n]$. \qed



\subsection{Sequential generation of $\G(n,\bfd)$}
\label{sec:d}

We define procedure \textsc{SeqSample-D}($\bfd$) which sequentially generates a random graph with distribution $\G(n,\bfd)$. This procedure is essentially the same as previously used in~\cite{gao2020sandwiching,gao2022sandwiching,klimovsova2023sandwiching,gao2020kim}. Let $\pr(jk\in \G(n,\bfd)\mid H)$ denotes the probability that $jk$ is an edge of $\G(n,\bfd)$ conditional on $H$ being a subgraph of $\G(n,\bfd)$.

\begin{algorithm}[H]
  \label{alg:d}
  \begin{algorithmic}[1]
    \Procedure{SeqSample-D}{$\bfd$}
    \State Let $G_0$ be the empty graph on $[n]$.
     \For{$i$ in $1, \dots,\norm{\bfd}_1/2$}
     \State Pick an edge $jk\in \binom{[n]}{2}\setminus G_{i-1}$ with probability proportional to $\pr(jk\in \G(n,\bfd)\mid G_{i-1})$.
     \State $G_i=G_{i-1}\cup \{jk\}$. 
     \EndFor
    \State Return $G_{\norm{\bfd}_1/2}$.
    \EndProcedure
    \end{algorithmic}
    \end{algorithm}

    Given $\bfd$, and an integer $0\le m\le \norm{\bfd}_1\jt{/2}$, let $\G(n,\bfd,m)$ denote a uniformly random subgraph of $\G(n,\bfd)$ with exactly $m$ edges. The following lemma follows by a simple counting argument and was proved in~\cite[Lemma 3]{gao2022sandwiching}. 
    
    \begin{lemma}\label{lem:G(n,d)}
    Let $0\le m \le \norm{\bfd}_1/2$ and $G_m$ be the graph obtained after $m$ iterations of \textsc{SeqSample-D}($\bfd$). Then $G_m\sim \G(n,\bfd,m)$.
    \end{lemma}
   
   Lemma~\ref{lem:G(n,d)} (with $m=\norm{\bfd}_1/2$) immediately implies the following corollary.
    
    \begin{corollary}\label{cor:G(n,d)}
Let $G$ be the output of $\textsc{SeqSample-D}({\bf d})$. Then $G \sim \G(n, {\bf d})$.
\end{corollary}

We need the following lemma, which follows by standard concentration results and can be found in~\cite[Lemma 4]{gao2022sandwiching}. \jt{We use $\Bin(K, p)$ to denote the Binomial distribution with $K$ trials and success probability $p$.} (Note that Lemma 4(c) of~\cite{gao2022sandwiching} only states the upper tail bound of part (c) below; but its proof gives both the upper and lower tail bounds.)

\begin{lemma}\label{lem:Bin-bounds} Let $Y \sim \Bin(K, p)$ for some positive integer $K$ and $p \in [0, 1]$. 
\begin{itemize}
    \item[(a)] For any $\eps\ge 0$, $\pr(|Y- p K| > \eps p K) \le 2e^{-\frac{\eps^2}{2+\eps} pK}$.
\item[(b)] If $p = j/K$ for some integer $j \in (0, K)$, then $\pr(Y = j) \ge \frac{1}{3} (p(1-p) K)^{-1/2}$.     
\item[(c)] Let $\I\sim \Po(\mu)$ for some $\mu>0$. Then, for any $\eps\ge 0$, $\pr(|\I -\mu |\ge \eps\mu)\le 2e^{-\frac{\eps^2}{2+\eps}\mu}$.
\end{itemize}
\end{lemma}

\jt{The following lemma is similar to~\cite[Lemma 6]{gao2022sandwiching}, whose proof easily extends to more general degree sequences $\bfd$ discussed in this paper through straightforward adjustments.} We include a short proof here. We will use this lemma to gain information on the remaining degree sequence of $\G(n,\bfd)$ when part of it has been constructed. \jt{Note that the probability lower bound $1-\exp\left(-\Omega(\xi^2 p_m d_j)+2\log n\right)$ below is not necessarily nonnegative. If it is negative then the assertion is trivially true.}

\begin{lemma}\label{lem:bounds}
Let $0\le m<\norm{\bfd}_1/2$ and $p_m = (\norm{\bfd}_1-2m)/\norm{\bfd}_1$. \jt{Let $d_j^{G_m}$ denote the degree of $j$ in $G_m$.} For any $\xi=\xi_n\in (0,1)$,
$|d_j-d^{G_m}_j-p_m d_j| \le \xi p_m d_j$ for all $j \in [n]$ with probability at least $1-\exp\left(-\Omega(\xi^2 p_m d_j)+2\log n\right)$.

\end{lemma}

\proof Take $G \sim \G(n, \bfd)$ and let $\textbf{h} = (h_1, \dots, h_n)$ be the degree sequence of the graph $H$ obtained by independently keeping every edge of $G$ with probability $p_m$.   Then, $h_j \sim \Bin(d_j, p_m)$ for every $j\in [n]$. By Lemma~\ref{lem:G(n,d)},  conditioned on the event that $|E(H)| = \norm{\bfd}_1/2-m$, $\textbf{h}$ has the same distribution as $\bfd-\bfd^{G_m}$. Therefore, by Lemma~\ref{lem:Bin-bounds} (a,b), for every $j\in [n]$, 
\begin{eqnarray*}
\pr(|d_j-d_j^{G_m}-p_m d_j| \ge \xi p_m d_j) &\le& \frac{\pr(|h_j-p_md_j| \ge \xi p_m d_j)}{\pr(|E(H)| = \norm{\bfd}_1/2-m)} \le \frac{2e^{-\frac{\xi^2}{2+\xi}p_md_j}}{\frac{1}{3}(p_m(1-p_m)\norm{\bfd}_1/2)^{-1/2}}\\
& \le& \sqrt{\norm{\bfd}_1} \cdot e^{-\Omega(\xi^2 p_m d_j)} = \exp\left(-\Omega(\xi^2 p_m d_j)+\log n\right),
\end{eqnarray*}
as $\norm{\bfd}_1\le n^2$. The lemma follows by taking the union bound over $j\in [n]$. \qed

\subsection{Couple the two sequential sampling procedures}
\label{sec:coupling}

Finally we couple the two aforementioned procedures and define procedure \textsc{Coupling}($\bfd$, $\lambda$, $\Lambda$) which sequentially constructs $\G(n,W)$ and $\G(n,\bfd)$ together.  \jt{The central idea of \textsc{Coupling} is that marginally, the constructions of $G_L\sim \G(n,W)$ and $G\sim \G(n,\bfd)$ follow precisely the procedures \textsc{SeqApprox-P} and \textsc{SeqSample-D} respectively. As in \textsc{SeqApprox-P}, edges $jk$ in $K_n$ are sequentially sampled independently with probability proportional to $d_jd_k$. If $jk$ was already added to $G$ then $jk$ is accepted by $G_L$  with probability  $\Lambda_{jk}$ as in \textsc{SeqApprox-P}. If $jk$ was not in $G$, then $jk$ is accepted to $G$ with some probability $\eta_{jk}$ so that the probability that $jk$ is selected and accepted by $G$ is proportional to the probability that $jk$ is in $\G(n,\bfd)$, conditional on the current construction of $G$, as desired in \textsc{SeqSample-D}. On the other hand, $jk$ is accepted by $G_L$ with probability $\Lambda_{jk}$ as required by \textsc{SeqApprox-P}. If $\eta_{jk}\ge \Lambda_{jk}$ then we can couple the two operations so that $jk$ is added to $G_L$ only when it is added to $G$. On the other hand, if $\eta_{jk}< \Lambda_{jk}$ then we simply sample $G_L\sim \G(n,W)$ and $G\sim \G(n,\bfd)$ independently. We can prove that the latter case happens rarely. The formal description of \textsc{Coupling} is given below.}

\begin{algorithm}[H]
  \label{alg:coupling}
  \begin{algorithmic}[1]
    \Procedure{Coupling}{$\bfd$, $\lambda$, $\Lambda$}
    \State Let $L_0$ and $G_0$ be empty graphs on vertex set $[n]$.
    \State Let $\I \sim \Po(\lambda)$.
    \For{$i$ in $1, \dots,\I$}
    \State Pick an edge $jk$ of $K_n$ with probability proportional to $d_jd_k$ (i.e.\ with probability $Q_{jk}$). 
    \If{$jk \in G_{i-1}$}
    \State $G_i = G_{i-1}$;
    \State $L_i = L_{i-1} \cup \{jk\}$ with probability $\Lambda_{jk}$,
    \State $L_i = L_{i-1}$ with probability $1-\Lambda_{jk}$.
    \Else 
    \State Let $\eta_{jk}^{(i)} = \frac{\rho_i(jk)}{\max_{h\ell\notin G_{i-1}}\rho_i(h\ell)}$, where $\rho_i(h\ell)=(d_hd_{\ell})^{-1}\pr(h\ell \in \G(n, {\bf d}) \mid G_{i-1})$. 
    \If{$\eta_{jk}^{(i)} <\Lambda_{jk}$}
    \State \textbf{Return} \textsc{IndSample}($\bfd,\lambda,\Lambda$)
    \Else
    \State $G_i = G_{i-1} \cup \{jk\}$ and $L_i = L_{i-1} \cup \{jk\}$ with probability $\Lambda_{jk}$,
    \State $G_i = G_{i-1} \cup \{jk\}$ and $L_i = L_{i-1}$ with probability $\eta_{jk}^{(i)}-\Lambda_{jk}$,
    \State $G_i = G_{i-1}$ and $L_i = L_{i-1}$ with probability $1-\eta_{jk}^{(i)}$.
    \EndIf
    \EndIf
    \EndFor
    \For{$i \ge \I+1$, while $G_{i-1}$ has fewer edges than $\G(n, {\bf d})$}
    \State Pick an edge $jk \not\in G_{i-1}$ with probability proportional to $\pr(jk \in \G(n, {\bf d}) \mid G_{i-1})$;
    \State $G_i = G_{i-1} \cup \{jk\}$.
    \EndFor
    \State \textbf{Return} $(G_L, G)$, where $G = G_i$ and $G_L= L_{\I}$. 
    \EndProcedure
    \item[]
    \Procedure{IndSample}{$\bfd,\lambda,\Lambda$}  
\State Independently sample $G\sim\G(n, \bfd)$, and $G_L\sim \G(n, f(\Lambda\odot Q))$ where $f(x)=1-\exp(-\lambda x)$. 
\State \textbf{Return} ($G_L$, $G$)
\EndProcedure
  \end{algorithmic}
\end{algorithm}

\begin{lemma}\label{lem:distribution} Let $(G_L,G)$ be the output of \textsc{Coupling}($\bfd,\lambda,\Lambda$).
\begin{enumerate}[(a)]
\item $G_L\sim \G(n, f(\Lambda\odot Q))$ where $Q=Q(\bfd)$ and $f(x)=1-\exp(-\lambda x)$. 
\item $G\sim \G(n,\bfd)$.
\item If \textsc{IndSample}($\bfd,\lambda,\Lambda$) is not called during the execution of \textsc{Coupling}($\bfd,\lambda,\Lambda$), then 
 $G_L\subseteq G$ in the output of \textsc{Coupling}($\bfd,\lambda,\Lambda$).
\end{enumerate}
\end{lemma}

\proof Parts (a,b) are trivially true if \textsc{IndSample}($\bfd,\zeta$) is called. Now assume that \textsc{IndSample}($\bfd,\zeta$) is not called. Part (c) follows directly by the coupling procedure, as every edge is added to $G_L$ only when it is, or has been added to $G$. For (a), note that if \textsc{IndSample}($\bfd,\zeta$) is not called then the edges are sequentially added to $G_L$ exactly as in \textsc{SeqApprox-P}($\bfd, \lambda, \Lambda$), and thus part (a) follows by Lemma~\ref{lem:P1}. For part (b), notice that in each step $i$, edge $jk$ is added to $G_{i-1}$ with probability $Q_{jk}\eta_{jk}^{(i)}$, which is proportional $\pr(jk\in \G(n,\bfd) \mid G_{i-1})$. Thus, the edges are added to $G$ exactly as in the execution of \textsc{SeqSample-D}($\bfd$), and part (b) follows by Corollary~\ref{cor:G(n,d)}. \qed 

\subsection{Specify $\lambda$ and $\Lambda$ for the coupling procedure}
\label{sec:parameters}
By the assumptions $J(\bfd)=o(\norm{\bfd}_1)$ and $\delta(\bfd)\gg \log n$, there exist $\xi,\zeta,\zeta'=o(1)$ that go to zero sufficiently slowly so that
\begin{align}
\zeta'&\gg \norm{\bfd}_1^{-1/3} \label{z'}\\
\zeta'\xi^2 & \gg \frac{\log n}{\delta(\bfd)} \label{z'-xi}\\
\zeta' & \gg \frac{J(\bfd)}{\norm{\bfd}_1}\label{z'2}\\
\zeta & \gg \frac{J(\bfd)}{\zeta'\norm{\bfd}_1}+\xi \label{z}.
\end{align}
\jt{To see why they exist, notice that there exists $\zeta'=o(1)$ satisfying~\eqn{z'2} since $J(\bfd)=o(\norm{\bfd}_1)$. Moreover, since $\delta(\bfd)\gg \log n$, we may assume that $\zeta'$ goes to 0 sufficiently slowly so that~\eqn{z'} is satisfied and there exists $\xi$ satisfying ~\eqn{z'-xi}. Finally, given~\eqn{z'2} and $\xi=o(1)$, the right hand side of~\eqn{z} is $o(1)$, and hence there exists $\zeta$ satisfying~\eqn{z}.}


Choose $\xi,\zeta,\zeta'$ that satisfy all the conditions above.
For the coupling procedure, we set
\begin{align}
\lambda &=(1-\zeta')\norm{\bfd}_1/2 \label{lambda}\\
\Lambda_{jk} &=(1-\zeta) \frac{\norm{\bfd}_1}{\norm{\bfd}_1+d_jd_k} \quad \mbox{for every $1\le j<k\le n$.} \label{Lambda}
\end{align}

\section{Proof of Theorem~\ref{thm:main2}}
\label{sec:proof}





Let $\Delta=\Delta(\bfd)$, $Q=Q(\bfd)$ and $P=P(\bfd)$.
By Lemma~\ref{lem:distribution},  it suffices to show that 
\begin{align}
 \pr(\textsc{IndSample}(\bfd,\zeta) \text{ is called}) &=o(1) \label{rejection} \\
\lambda \Lambda\odot Q &= (1+o(1))P, \label{Q}
\end{align}
as $1-\exp(-(1+o(1))P_{ij})=(1+o(1))(1-e^{-P_{ij}})$ for every $ij$. \medskip

\noindent{\em Proof of~\eqn{Q}.} For every $1\le j<k\le n$, 
\[
Q_{jk}= \frac{d_jd_k}{\sum_{1\le h<\ell\le n} d_hd_{\ell}}= \frac{d_jd_k}{\frac{1}{2}(\norm{\bfd}_1^2-\sum_{i= 1}^n d_i^2)} = \frac{2d_jd_k}{\norm{\bfd}_1^2- O(\Delta) \norm{\bfd}_1}.  
\]
Since $\lambda = (1-\zeta')\norm{\bfd}_1/2$ by~\eqn{lambda}, and by~\eqn{Lambda}
\[
\lambda \Lambda_{jk} Q_{jk} = \frac{(1-\zeta') \norm{\bfd}_1}{2}  \frac{(1-\zeta)\norm{\bfd}_1}{\norm{\bfd}_1+d_jd_k}  \frac{2d_jd_k}{\norm{\bfd}_1^2(1+ O(\Delta/\norm{\bfd}_1) } =\left(1+O\left(\zeta+\zeta'+\frac{\Delta}{\norm{\bfd}_1}\right)\right) \frac{d_jd_k}{\norm{\bfd}_1+d_jd_k}, 
\]
and~\eqn{Q} follows by noting that $\zeta,\zeta',\Delta/\norm{\bfd}_1=o(1)$. \medskip


\noindent {\em Proof of~\eqn{rejection}.} By Lemma~\ref{lem:Bin-bounds}(c) \jt{(with $\mu=\lambda$ and $\eps=\lambda^{-1/3}$)} and~\eqn{lambda}, a.a.s.\ $\I=(1+O(\lambda^{-1/3}))\lambda = (1-\zeta'+O(\norm{\bfd}_1^{-1/3}))\norm{\bfd}_1/2$. It suffices then to show that a.a.s.\ throughout the execution of  \textsc{Coupling}($\bfd,\lambda,\Lambda$), $\eta^{(i)}_{x_i}\ge \Lambda_{x_i}$ for every $1\le i\le \I$, where $x_i$ is the random edge sampled in the $i$th iteration. Let $m_i$ be the number of edges in $G_i$. Since a.a.s.\ $\I=(1-\zeta'+O(\norm{\bfd}_1^{-1/3}))\norm{\bfd}_1/2$ and $m_i\le \I$ for every $1\le i\le \I$, it follows then by~\eqn{z'} that a.a.s.\ 
\begin{equation}
p_{m_i}\ge p_{m_\I}\ge \zeta'/2 \quad \mbox{for every $1\le i\le \I$}, \label{p_m}
\end{equation}
 where $p_{m_i}$ is defined by $1-2m_i/\norm{\bfd}_1$ as in Lemma~\ref{lem:bounds}. 
 
 By Lemma~\ref{lem:bounds} (with $\xi$ chosen in Section~\ref{sec:parameters}),~\eqn{p_m} and~\eqn{z'-xi}, \jt{with probability at least 
 $$
 1-\exp\left(-\Omega(\xi^2 p_m d_j)+2\log n\right)\ge  1-\exp\left(-\Omega(\xi^2 \zeta' \delta(\bfd))+2\log n\right) =  1-\exp\left(-\omega(\log n\right), 
 $$
 we have that $|d_j-d^{G_i}_j-p_{m_i} d_j| \le \xi p_{m_i} d_j$ for all $j \in [n]$.}
Now take the union bound over all the $\I\le \norm{\bfd}_1/2$ steps, we conclude that 
\begin{equation}
d_j - d_j^{G_i}=(1+O(\xi)) p_{m_i} d_j\quad \mbox{for every $j\in [n]$, and for every $i\le \I$.} \label{eq:residual} 
\end{equation}

Next, we estimate $\pr(jk\in \G(n,\bfd) \mid G_i)$. \jt{Note that Proposition~\ref{p:edgeprob} gives $\pr(jk\in \G(n,\bfd) \mid G_i)$ when $G_i=\emptyset$. Here, we apply another corollary of~\cite[Theorem 1]{gao2020subgraph}, given in Theorem~\ref{lem:edgeprob} below, which
 estimates the edge probabilities in $\G(n,\bfd)$ when conditioned on a set of edges $H$ being present in $\G(n,\bfd)$. }

For two degree sequences ${\bfd}$ and $\bfg$, we say $\bfd\preceq \bfg$ if $d_i\le g_i$ for every $i\in [n]$. Given a graph $H$ on $[n]$, let $\bfd^H$ denote the degree sequence of $H$.

\begin{theorem}\label{lem:edgeprob}
Suppose $H$ is a graph on $[n]$ with ${\bf d}^{H}\preceq {\bf d}$, and let ${\bf t} = {\bf d}-{\bf d}^H$. Suppose that $J(\bfd)=o(\norm{\bft}_1)$ and
suppose $jk \not\in H$. Then
\[
\pr(jk \in \G(n, {\bf d}) \mid H) =\left(1+O\left(\frac{J(\bfd)}{\norm{\bft}_1}\right)\right) \frac{t_jt_k}{\norm{\bft}_1+t_jt_k},
\]
where $\pr(jk \in \G(n, {\bf d}) \mid H)$ denotes the probability that $jk\in \G(n,\bfd)$ conditional on the event that $\G(n,\bfd)$ contains $H$ as a subgraph.
\end{theorem}

\jt{We first verify that the assumption $J(\bfd)=o(\norm{\bft}_1)$ of Theorem~\ref{lem:edgeprob} is satisfied for every $i\le\I$ with $H=G_i$. By~\eqn{p_m} we may assume that $p_{m_i}\ge \zeta'/2$ for every $1\le i\le \I$ (which holds a.a.s.). Let $\bft_i=\bfd-\bfd^{G_i}$. Then, 
\begin{equation}
\norm{\bft_i}_1=\norm{\bfd}_1-2m_i=\norm{\bfd}_1(1-2m_i/\norm{\bfd}_1)=\norm{\bfd}_1p_{m_i}\ge \norm{\bfd}_1\zeta'/2.\label{eq:ti}
\end{equation}
By~\eqn{z'2}, $J(\bfd)\ll \zeta' \norm{\bfd}_1$ and consequently $J(\bfd)/\norm{\bft_i}_1=O(J(\bfd)/\zeta'\norm{\bfd}_1)=o(1)$ for every $1\le i\le \I$.}
  By Theorem~\ref{lem:edgeprob},~\eqn{eq:residual} and~\eqn{eq:ti}, a.a.s.\ 
\[
\pr(jk\in \G(n,\bfd) \mid G_i)=\left(1+O\left(\frac{J(\bfd)}{p_{m_i}\norm{\bfd}_1} + \xi\right)\right) \frac{ p_{m_i} d_j d_k }{ \norm{\bfd}_1+p_{m_i} d_jd_k },\quad \mbox{ for every $1\le i\le \I$,}
\]
 where $J(\bfd)/p_{m_i}\norm{\bfd}_1 + \xi = O(J(\bfd)/\zeta'\norm{\bfd}_1 + \xi)=o(1)$ as shown before. 
 Recall that 
 \[
 \rho_i(jk)=\frac{\pr(jk\in \G(n,\bfd) \mid G_i)}{d_jd_k}.
  \]
  It follows that
  \[
   \rho_i(jk)=\left(1+O\left(\frac{J(\bfd)}{\zeta'\norm{\bfd}_1} + \xi\right)\right) \frac{ p_{m_i}}{ \norm{\bfd}_1+p_{m_i} d_jd_k }.
  \]
 Notice that
 \[
 \frac{ p_{m_i}}{ \norm{\bfd}_1+p_{m_i} d_jd_k } \le \frac{ p_{m_i}}{ \norm{\bfd}_1} \quad\mbox{for all $jk\notin G_{i-1}$}.
 \]
 Thus, 
 $$
 \max_{h\ell\notin G_{i-1}}\rho_i(h\ell) \le \left(1+O\left(\frac{J(\bfd)}{\zeta'\norm{\bfd}_1} + \xi\right)\right) p_{m_i}/\norm{\bfd}_1.
 $$
 Hence, 
 \[
 \eta^{(i)}_{jk}\ge \left(1+O\left(\frac{J(\bfd)}{\zeta'\norm{\bfd}_1} + \xi\right)\right) \frac{\norm{\bfd}_1}{\norm{\bfd}_1+p_{m_i}d_jd_k}\ge (1-\zeta)\frac{\norm{\bfd}_1}{\norm{\bfd}_1+d_jd_k},
 \]
 by~\eqn{z}. Now~\eqn{rejection} follows. \qed \medskip

\noindent {\bf Acknowledgement} We thank two anonymous referees for their valuable advice in improving the presentation of the paper.


\fix{
\section*{Appendix}

We prove a few properties of degree sequences $\bfd$ where $J(\bfd)=o(\norm{\bfd}_1)$. 
\begin{proposition}\label{p:properties}
Let $\bfd=\bfd_n$ be a sequence of degree sequences such that $J(\bfd)=o(\norm{\bfd}_1)$. Let $\Delta=\max\{d_i,i\in[n]\}$. Then,
\begin{enumerate}[(a)]
\item $\Delta=o(n)$.
\item For every $\eps>0$, $ |\{ij: d_id_j\ge \eps \norm{\bfd}_1\}|\le \eps n^2$ for all sufficiently large $n$.
\end{enumerate}
\end{proposition}

\proof
 For (a), suppose on the contrary that $\Delta\ge \delta n$ for some absolute constant $\delta>0$. Then, 
 \[
 \norm{\bfd}_1=\sum_{i=1}^n d_i \le \delta^{-1} \sum_{i=1}^{\delta n} d_i \le \delta^{-1} J(\bfd),
 \]
 contradicting with $J(\bfd)=o(\norm{\bfd}_1)$.
 
 For (b), note that $J(\bfd)\ge \Delta d_{\Delta}$ and thus $\Delta d_{\Delta} = o(\norm{\bfd}_1)$. It follows that for every $i,j \in [n]$ such that $j\ge \Delta$, $d_id_j \le \Delta d_{\Delta} = o(\norm{\bfd}_1)$. Thus, part (b) follows by $\Delta=o(n)$ from part (a). \qed \medskip
 
 \begin{example}
We construct an example of degree sequences $\bfd\preceq \bfg$ for which there exists $jk$ such that $\pr(jk\in \G(n,\bfd))= 1-o(1)$ and $\pr(jk\in \G(n,\bfg))=o(1)$. 
Let $\bft$ be a degree sequence where $t_1=t_2=n^{2/3}$, and $t_3=\cdots=t_n=1$ (without loss of generality we assume that $n^{2/3}$ is an integer and $\norm{\bft}_1$ is even). Let $\bfg=(n-1){\bf 1}-\bft$, and let $\bfd$ be $(n-2n^{2/3})$-regular. Obviously, $\bfd\preceq \bfg$. By symmetry, the probability that vertices 1 and 2 are adjacent in $\G(n,\bfd)$ is equal to $(n-2n^{2/3})/(n-1)=1-o(1)$. By Proposition~\ref{p:edgeprob}, the probability that these two vertices are adjacent in $\G(n,\bft)$ is $1-o(1)$, and consequently, their adjacency probability in $\G(n,\bfg)$ is $o(1)$, as $\G(n,\bfg)$ has the same distribution as the complement of $\G(n,\bft)$.
 \end{example}
 
 }

\end{document}